\documentclass[a4paper]{jpconf}
\usepackage{graphicx}
\begin{document}
\title{ On a family of complex algebraic surfaces
of degree 3n }
 \author{Juan Garc\'{\i}a Escudero}

\address{Universidad de Oviedo. Facultad de Ciencias,
 33007 Oviedo, Spain}

\begin{abstract}
We study a class of algebraic surfaces of degree $3n$ in the complex projective space with only ordinary double points. They are obtained by  using bivariate polynomials with complex coefficients related to the generalized cosine associated to the affine Weyl group of the root system $A_{2}$.
\end{abstract}
%

  \bigskip\par  

In \cite{chm92}  Chmutov introduced surfaces $V_{d}$ of degree $d$ in the complex projective space ${\bf{P}}^{3}({\bf{C}})$ given by the affine equations $P_{d}(u,v)+T_{d}(w)=0$, where $P_{d}$ are the generalized Chebyshev polynomials or folding polynomials \cite{hof88,wit88} associated to the affine Weyl group $\tilde{W}(A_{2})$ of the root system $A_{2}$, and $T_{d}$ is a Chebyshev polynomial. In this work we study a class of surfaces with degree $3n$ which have more singularities than the Chmutov surfaces of the same degree.
\par
For $d=3n$, $n=1,2,3,...$ the surfaces $U_{d}$ are given by the affine equations
  \begin{equation}
Q_{d}(u,v)+T_{d}(w)=0
\end{equation}
\par\noindent
where $T_{d}$ is the Chebyshev polynomial with two critical values 2 and 3 and $Q_{d}$ is defined as follows.
  The generalized cosine associated to the affine Weil group $\tilde{W}(A_{2})$ is $h(u,v)=(h_{1},h_{2})$, where 
  \begin{equation}
h_{1}(u,v):=e^{-2\pi iu}+e^{-2\pi i v}+e^{2\pi i(u+v)}, h_{2}:=\bar{h}_{1}
\end{equation}
 The polynomials appearing in the surfaces $V_{d}$ are $P_{d}(u,v):=C_{d} \circ h^{-1}(u,v)$, where $C_{d}(u,v):=h_{1,d}(u,v)+h_{2,d}(u,v), h_{1,d}(u,v)=h_{1}(du,dv), h_{2,d}=\bar{h}_{1,d}$. The number of non-degenerate singularities of $V_{d}$ for $d=3n$ is \cite{chm92}
   \begin{equation}
{d \choose 2}\lfloor \frac{d}{2}\rfloor+(\frac{d^{2}}{3}-d)\lfloor \frac{d-1}{2}\rfloor
\end{equation}
 \par
 For the construction of the surfaces $U_{d}$ we use 
  \begin{equation}
g_{1,d}(u,v):=e^{-2\pi i(du+\frac{2}{3})}+e^{-2\pi i(dv-\frac{1}{3})}+e^{2\pi i(d(u+v)+\frac{1}{3})}, g_{2,d}:=\bar{g}_{1,d}
\end{equation}
\par\noindent
  The bivariate polynomials in eq(1) are defined by $Q_{d}(u,v):=H_{d} \circ h^{-1}(u,v)$, where
     \begin{equation}
  H_{d}(u,v):=g_{1,d}(u,v)+g_{2,d}(u,v)= 2 {\rm cos}(2\pi d u-\frac{2\pi }{3})+2 {\rm cos}(2\pi d v-\frac{2\pi }{3})+2 {\rm cos}(2\pi d (u+v)+\frac{2\pi }{3})
  \end{equation}
   \bigskip\par  
  {\it{  Lemma.}} 
 The polynomial $Q_{d}$ has $\frac{d(d-3)}{6}$ critical points with critical value 6, ${d \choose 2}$ critical points with critical value -2 and $\frac{d^{2}}{3}-d+1$ critical points with critical value -3. All the critical points of $Q_{d}$ are non-degenerate. 
   \bigskip\par  
  {\it{  Proof.}} 
   We restrict $H_{d}$ and $h$ on the plane with real coordinates. The Jacobi matrices satisfy $J(H_{d})=J(Q_{d}) \circ J(h)$. The zeros of  the Jacobian determinant of $h$ are the sides of the triangle $\Delta$ whose interior is given by $u-v>0,u+2v>0,2u+v<1$, which is the fundamental region of $\tilde{W}(A_{2})$. All the points from one orbit of $\tilde{W}(A_{2})$ are mapped by $h$ into a single point, and the images of the critical points of $H_{d}$ from the interior of $\Delta$ by $h$ are the critical points of $Q_{d}$.
    \par
   A direct computation of the critical points of $H_{d}$ leads to three cases. In the following list we indicate the critical value $\zeta$ and the number of points $N_{\zeta}$  corresponding to $\zeta$  inside $\Delta$ by $(\zeta,N_{\zeta})$:
  \bigskip\par 
   \par\noindent
 a) (6,$\frac{d(d-3)}{6}$); $u=\frac{3k-2}{d},v=\frac{3l+1}{d}$, for $k,l \in \bf{Z}$.
 \bigskip\par 
 \par\noindent
 b) (-3,$\frac{d^{2}}{3}-d+1$); $u=\frac{k-2}{3d},v=\frac{l+1}{3d}$, with $k=3m+1,l=3n+1$ or $k=3m+2,l=3n+2$, $m,n \in \bf{Z}$.
 \bigskip\par 
  b1) (-3,$1+\frac{d(d-3)}{6}$); $u=\frac{3m-1}{3d},v=\frac{3n+2}{3d}$, for $m,n \in \bf{Z}$.     
 \bigskip\par 
  b2) (-3,$\frac{d(d-3)}{6}$); $u=\frac{m}{d},v=\frac{n+1}{d}$, for $m,n \in \bf{Z}$. 
   \bigskip\par 
 \par\noindent
 c) (-2,${d \choose 2}$); $u=\frac{3k-4}{6d},v=\frac{3l+2}{6d}$ with $k$ or $l$ odd.
 \bigskip\par 
  c1) (-2,$\frac{d(d-1)}{3}$); $u=\frac{6m-1}{6d},v=\frac{3n+2}{6d}$, for $m,n \in \bf{Z}$.
 \bigskip\par 
  c2) (-2,$\frac{d(d-1)}{6}$); $u=\frac{6m-4}{6d},v=\frac{6n+5}{6d}$, for $m,n \in \bf{Z}$.
   \bigskip\par

The hessian matrix of $H_{d}$
        \par
\[ \left(
    \begin{array}{ll}
   -8\pi^2  d^2({\rm cos}(2\pi d u-\frac{2\pi }{3})-{\rm cos}(2\pi d (u+v)+\frac{2\pi }{3})) & -8\pi^2  d^2 {\rm cos}(2\pi d (u+v)+\frac{2\pi }{3}) \\
-8\pi^2  d^2 {\rm cos}(2\pi d (u+v)+\frac{2\pi }{3})  &  -8\pi^2  d^2({\rm cos}(2\pi d v-\frac{2\pi }{3})-{\rm cos}(2\pi d (u+v)+\frac{2\pi }{3})) \\

    \end{array}
    \right) \]
  \par\noindent
has full rank in all the critical points, hence they are non-degenerate.    
\bigskip\par  
In Figs.1,2 we can see the critical points of $H_{d}$ inside $\Delta$ for $d=6,9$. Critical points with critical values 6,-2,-3 are represented by $\circ, \ast, \bullet$ respectively. The distance between two consecutive lines in the $(u,v)$ oblique coordinate system is $\frac{1}{6d}$.
\bigskip\par 
   {\it{  Theorem.}}   The number of singular points of $U_{d}$ is  ${d \choose 2}\lfloor \frac{d}{2}\rfloor+(\frac{d^{2}}{3}-d+1)\lfloor \frac{d-1}{2}\rfloor$. All the singular points are non-degenerate.
 The surface cannot have singular points at infinity.
       \bigskip\par 
  {\it{  Proof.}}  The Chebyshev polynomials $T_{d}(w)$ have $\lfloor \frac{d}{2}\rfloor$ critical points with critical value 2 and $\lfloor \frac{d-1}{2}\rfloor$ critical points with critical value 3. The surface is singular at the points where the sum of the critical values of $T_{d}(w)$ and $Q_{d}(u,v)$ is zero. The result for the number of non-degenerate singularities of $U_{d}$ follows then from the Lemma: $Q_{d}$ has ${d \choose 2}$ critical points with critical value -2 and $\frac{d^{2}}{3}-d+1$ critical points with critical value -3.
  \par
  In the Lemma we have also shown that the number of distinct critical points of $Q_{d}$ is $(d-1)^2$, therefore $U_{d}$ can not have singular points at infinity.  
    \bigskip\par 
    Consider a surface of degree $d$ in ${\bf{P}}^{3}({\bf{C}})$ with $N(d)$ double points and no other singularities, and let  $\mu(d)={\rm max}$ $N(d)$. Then we have 
    \bigskip\par 
    {\it{Corollary}}.  $\mu(3n) \ge {3n \choose 2}\lfloor \frac{3n}{2}\rfloor+(3n^{2}-3n+1)\lfloor \frac{3n-1}{2}\rfloor$.
       \bigskip\par 
    We notice that $U_{3n}$ has $\lfloor \frac{3n-1}{2}\rfloor$ more singularities than $V_{3n}$. Also of interest are hypersurfaces in ${\bf{P}}^{4}({\bf{C}})$ with affine equations.
    
      \begin{equation}
Q_{3n}(u_{1},u_{2})-Q_{3n}(u_{3},u_{4})=0
\end{equation}
  They have  $(\frac{3n(3n-1)}{2})^{2}+(3n(n-1)+1)^{2}+(\frac{3n(n-1)}{2})^{2}$ non-degenerate singularities. We find $3n(n-1)$ more singularities than in the Chmutov hypersurfaces $P_{3n}(u_{1},u_{2})-P_{3n}(u_{3},u_{4})=0$.
  \bigskip\par 
Real variants of $V_{d}$ were studied in \cite{bre08}, and the authors showed that the known lower bounds for the maximum number of ordinary double points on a surface of degree $d$ can be realized with only real singularities. Recently we have shown that a construction connected to the derivation of substitution tilings \cite{esc08} can be used for the generation of algebraic surfaces with many real nodes \cite{esc11}. One of the two types of surfaces obtained is equivalent to real variants of $V_{d}$. The other type consists in surfaces of degree $3n$ which have the same number of singularities as the surfaces presented in this work. In fact they are related to the real variants of  $U_{3n}(u,v,w)$ with $u=x+iy,v=\bar{u},w=z$, for $x,y,z \in {\bf{R}}$.
\bigskip\par

{ \bf{References}}
 \bigskip\par

\newpage

 \begin{figure}[h]
 \includegraphics[width=20pc]{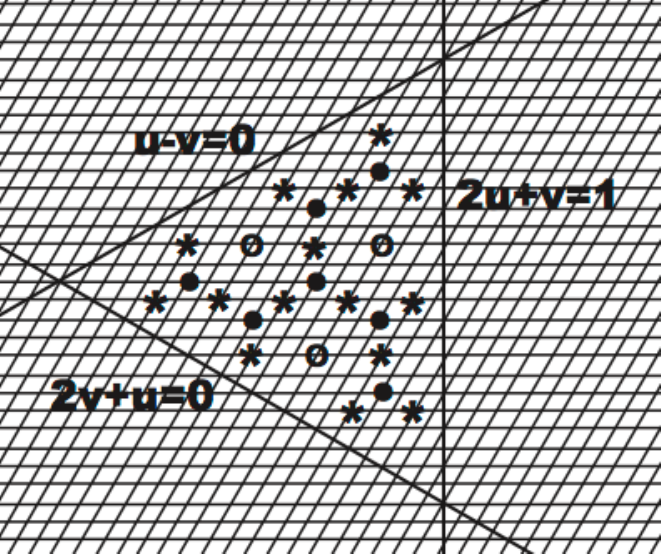}
\caption{\label{label} The critical points of $H_{6}$ inside the fundamental region $\Delta$ of the affine Weyl group $\tilde{W}(A_{2})$.}
\end{figure}

 \begin{figure}[h]
 \includegraphics[width=20pc]{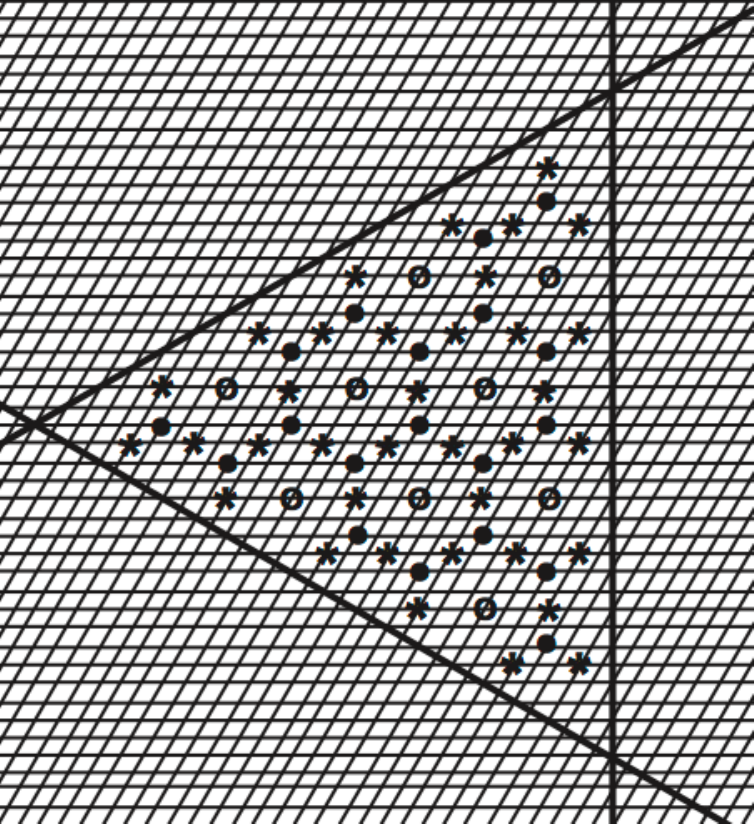}
\caption{\label{label} The critical points of $H_{9}$ inside $\Delta$.}
\end{figure}


\begin{thebibliography}{9}



\bibitem{bre08} S.Breske, O.Labs, D. van Straten, {\em Real line arrangements and surfaces with many real nodes.} In Geometric modeling and algebraic geometry.  Springer. Berlin. (2008) 47-54.

\bibitem{chm92}  S.V.Chmutov, {\em Examples of projective surfaces with many singularities.} J.Algebr.Geom.  {\bf 1} (1992) 191-196.

\bibitem{esc08} J.G. Escudero,  {\em Random tilings of spherical 3-manifolds. } J.Geom.Phys. {\bf 58} (2008)  1451-1464.

\bibitem{esc11} J.G.Escudero, {\em A construction of algebraic surfaces with many real nodes.} http://arxiv.org/abs/1107.3401 (2011).

  \bibitem{hof88}
  M.E. Hoffman, D. Withers, {\em Generalized Chebyshev polynomials associated with affine Weyl groups.}  Trans. Amer. Math.Soc.{\bf 282} (1988)  555-575.
  
    \bibitem{wit88}
  D. Withers, {\em Folding polynomials and their dynamics.}  Amer.Math.Monthly. {\bf 95} (1988)  399-413.



\end{thebibliography}
\end{document}